\documentclass[reqno]{amsart}
\usepackage{amscd}
\usepackage{amsmath}
\usepackage{amsfonts}
\usepackage{amssymb}

\theoremstyle{plain}

\newtheorem{proposition}{Proposition}

\numberwithin{equation}{section}
\newcommand{\R}{\mathbb{R}}
\newcommand{\C}{\mathbb{C}}

\begin{document}
\title[Nonlinear Filtration Equation]{Global Actions of the Lie Symmetries of the Nonlinear Filtration Equation}
\author[J. Franco]{Jose A. Franco}
\address{Department of Mathematics and Statistics \\ University of North Florida \\ 1 UNF Drive \\ Jacksonville, FL 32224}
\email{jose.franco@unf.edu}

\subjclass{22E25, 22E27, 22E70}
\keywords{Globalizations, solvable Lie groups, nonlinear filtration equation, parabolic induction}

\begin{abstract}
The classification of the Lie point symmetries of the nonlinear filtration equation gives the generic case and three special cases. By restricting to a special class of functions, we show that the Lie symmetries of the nonlinear filtration equation exponentiate to a global action of a solvable Lie group in the generic case and two of the three special cases. We show that the action of the Lie point symmetries cannot be globalized for the third special case.  
\end{abstract}
\maketitle

\section{Introduction}
The use of Lie groups in the analysis of differential equations is a powerful tool that has been vastly studied in mathematics (c.f. \cite{ CNC, bluman2002symmetry, CRC2,Olver} and the references therein to mention just a few). Algorithms like the prolongation algorithm provide infinitesimal generators whose actions exponentiate to one-parameter local Lie groups (c.f. \cite{Olver}). However, it is well-known that these actions generally do not exponentiate to a global action of a Lie group. Therefore the tools from representation theory are not readily available to use as these apply mostly to global Lie groups. Nonetheless, M. Craddock discovered that by restricting to a special subset of the solution space of a partial differential equation, the local actions generated by its symmetry Lie algebra could be globalized (c.f. \cite{Craddock}). In \cite{Franco2, Franco3,Sepanski3, Sepanski} the globalization problem is solved by using the representation theoretic method of parabolic induction. Using this method spaces of solutions to the different partial differential equations are constructed. These spaces consist of smooth functions on which the actions of the Lie symmetries globalize. In all these cases the PDE's were linear. However, these constructions can be used in some nonlinear equations as well. Some examples can be found in \cite{Franco, Sepanski4, Sepanski2}.

In this article we use induced representations of solvable groups to study the globalization problem for the the nonlinear filtration equation
\begin{equation}\label{NLFE}
v_t=k(v_x)v_{xx}
\end{equation}
with $k'(v_x)\neq 0$. The solution to its classification problem is well-known and readily available in \cite{CNC,Akhatov}. For the general case, the symmetry Lie algebra $\mathfrak{g}$ is spanned by the following infinitesimal generators:
\begin{align}
X_1&=\frac{\partial}{\partial t}, & X_2&=\frac{\partial}{\partial x}, &X_3&=\frac{\partial}{\partial v}, & X_4&=2t\frac{\partial}{\partial t}+x\frac{\partial}{\partial x}+v\frac{\partial}{\partial v}.
\end{align}
The algebra extends for $k(v_x)=e^{v_x}$, $k(v_x)=(v_x)^n$, and $k(v_x)=\frac{e^{n \arctan(v_x)}}{1+v_x^2}$.
The main goal of this article is to study the globalization problem for the local actions of the one-parameter groups generated by $\mathfrak{g}$ and its extensions. For the generic case and for   $k(v_x)=e^{v_x}$ and $k(v_x)=(v_x)^n$ we successfully construct a space of functions on which the action of the symmetry group corresponding to each equation globalizes.  When $k(v_x)=\frac{e^{n \arctan(v_x)}}{1+v_x^2}$, however, we show that this is not possible to globalize the action of the one parameter subgroups on a non-trivial set of functions.
%----------------------------------------------- Section 2 -------------------------------------------------------------------
\section{Generic Case}\label{genericcase}
Consider the solvable group $$G_1=\left\{\begin{pmatrix} q^2 & 0 & 0 & t \\ 0 & q &0 & x \\  0 & 0 &q & s \\ 0 & 0 &0 & 1 \end{pmatrix} \ \bigg| \ s,t,x\in \R \text{ and } q \in \R^{>0} \right\}$$
and the subgroup
$$P_1=\left\{\begin{pmatrix} q^2 & 0 & 0 & 0 \\ 0 & q &0 & 0 \\  0 & 0 &q & s \\ 0 & 0 &0 & 1 \end{pmatrix} \ \bigg| \ s\in \R \text{ and } q \in \R^{>0} \right\}.$$
Define a character $\chi_{1}:P_1\to \C^\times$ by $$\chi_{1}\begin{pmatrix} q^2 & 0 & 0 & 0 \\ 0 & q &0 & 0 \\  0 & 0 &q & s \\ 0 & 0 &0 & 1 \end{pmatrix}=q.$$
The induced representation associated to this character would be
$$I_{1}:=\text{Ind}_{P_1}^{G_1}(\chi_{1}):=\left\{\varphi \in C^\infty(G) \ | \  \varphi(gp)=\chi_{1}(p)^{-1} \varphi(g) \text{ for all } g\in G_1 \text{ and } p\in P_1\right\}.$$
This space carries the structure of a $G_1$-representation, with the action of the group given by:
$$g_1\cdot \varphi(g_2)=\varphi(g_1^{-1}g_2).$$
For the rest of this section let 
\begin{align}\label{GroupElements} g_1=&\begin{pmatrix} q_1^2 & 0 & 0 & t_1 \\ 0 & q_1 &0 & x_1 \\  0 & 0 &q_1 & s_1 \\ 0 & 0 &0 & 1 \end{pmatrix}& \text{ and } & & g_2=&\begin{pmatrix} q_2^2 & 0 & 0 & t_2 \\ 0 & q_2 &0 & x_2 \\  0 & 0 &q_2 & s_2 \\ 0 & 0 &0 & 1 \end{pmatrix}. \end{align} 
Notice that the subgroup $$N_1:=\left\{n_{t,x}:=\begin{pmatrix} 1 & 0 & 0 &t  \\ 0 & 1&0 & x \\  0 & 0 &1 & 0 \\ 0 & 0 &0 & 1 \end{pmatrix} \ \bigg| \ t,x\in \R \right\}$$ is isomorphic to   $\R^2$. Using this isomorphism, we can realize $\varphi \in I_{1}$ as a smooth map on $\R^2$ by restriction of domain. That is, for $\varphi \in I_{1}$ define $f\in C^\infty(\R^2)$ by $f(t,x)=\varphi(n_{t,x}).$
\begin{proposition}Let
$$I'_{1}:=\left\{ f\in C^\infty(\R^2) \ | \ f(t,x)=\varphi(n_{t,x}) \text{ for some } \varphi \in I_1  \text{ and for all } t,x\in \R^2 \right\}.$$
Then, $I_{1}\cong I'_{1}$ as $G_1$-representations. 
\begin{proof}
Notice that we can factor
$$
\begin{pmatrix} q^2 & 0 & 0 & t \\ 0 & q &0 & x \\  0 & 0 &q & s \\ 0 & 0 &0 & 1 \end{pmatrix}=\begin{pmatrix} 1 & 0 & 0 &t  \\ 0 & 1&0 & x \\  0 & 0 &1 & 0 \\ 0 & 0 &0 & 1 \end{pmatrix}\begin{pmatrix} q^2 & 0 & 0 & 0 \\ 0 & q &0 & 0 \\  0 & 0 &q & s \\ 0 & 0 &0 & 1 \end{pmatrix} .
$$
Thus, by the definition of $I_1$ and $\chi_1$,
$$\varphi\begin{pmatrix} q^2 & 0 & 0 & t \\ 0 & q &0 & x \\  0 & 0 &q & s \\ 0 & 0 &0 & 1 \end{pmatrix}=q^{-1}\varphi\begin{pmatrix} 1 & 0 & 0 &t  \\ 0 & 1&0 & x \\  0 & 0 &1 & 0 \\ 0 & 0 &0 & 1 \end{pmatrix}=q^{-1}f(t,x).$$
So, the restriction map is a bijection and clearly linear, hence an isomorphism of vector spaces. 
Then, we define the action of $G_1$ on $I'_1$ by
\begin{equation}\label{LinearAction}g_1\cdot f(t,x)=q_1f\left(\frac{t-t_1}{q_1^2},\frac{x-x_1}{q_1} \right).\end{equation}

It is easy to check that this action makes the restriction map an intertwining map completing the proof.
\end{proof}
\end{proposition} 
An extra piece of information that we can obtain from the proof of this proposition is that given an arbitrary smooth function $f(t,x)$ we can construct a corresponding $\varphi \in I_1$. Therefore, $I_{1}\cong C^\infty(\R^2)$.

Notice the assignment
\begin{align}\label{LAgenerators}
X_1 &\leftrightarrow \xi_1:= \begin{pmatrix}0 & 0 & 0 & 1 \\ 0 & 0 &0 & 0 \\  0 & 0 &0 & 0 \\ 0 & 0 &0 & 0 \end{pmatrix}&   X_2 \leftrightarrow&\xi_2:= \begin{pmatrix}0 & 0 & 0 & 0 \\ 0 & 0 &0 & 1 \\  0 & 0 &0 & 0 \\ 0 & 0 &0 & 0 \end{pmatrix} \nonumber \\
X_3 &\leftrightarrow\xi_3:=  \begin{pmatrix}0 & 0 & 0 & 0 \\ 0 & 0 &0 & 0 \\  0 & 0 &0 & 1 \\ 0 & 0 &0 & 0 \end{pmatrix}&   X_4\leftrightarrow  &\xi_4:=\begin{pmatrix}2 & 0 & 0 & 0 \\ 0 & 1 &0 & 0 \\  0 & 0 &1 & 0 \\ 0 & 0 &0 & 0 \end{pmatrix}
\end{align}
gives a Lie algebra isomorphism between $\mathfrak{g}$ and $\text{Lie}(G_1)$. Moreover, differentiating the action of the one-parameter subgroups generated by these elements of the Lie algebra according to \eqref{LinearAction}, gives the actions of $X_1, X_2,$ and $X_4$ on $I'_1$. However, the one parameter group generated by $\xi_3$ is the trivial group containing only the identity, but $X_3$ generates a one parameter subgroup isomorphic to $\R$. Therefore, we need to modify the  action of $G_1$ in order to obtain a globalization of the local action of the one parameter groups generated by $X_1,X_2,X_3,$ and $X_4$. To do this, we will let $\theta_1:G_1\times I'_{1}\to I'_{1}$ be defined by $$\theta_1(g_1)\cdot f(t,x)= f(t,x)+s_1,$$ with $g_1$ defined by \eqref{GroupElements}. It is straightforward to see that 
\begin{align*}
\theta_1(g_1)\cdot (\theta_1(g_2)\cdot f) (t,x)& = f(t,x)+s_1+s_2\\
(\theta_1(g_1 g_2)\cdot f)(t,x)&=f(t,x)+s_1+q_1 s_2.
\end{align*}
So, $\theta_1$ does not define an action of $G_1$ on $I'_{1}$. However, we can still use $\theta_1$ to define an action of $G_1$ on $I'_{1}$ that globalizes the action of $\mathfrak g$ on $I'_{1}$. To do this, we notice that
\begin{align}\label{ThetaProduct}\theta_1(g_1) \cdot (g_1\cdot (\theta_1(g_2)\cdot(g_1^{-1}\cdot f)))(t,x) &=g_1\cdot (\theta_1(g_2)\cdot(g_1^{-1}\cdot f))(t,x)+s_1 \nonumber
\\ \quad &=q_1(\theta_1(g_2)\cdot(g_1^{-1}\cdot f))\left(\frac{t-t_1}{q_1^2},\frac{x-x_1}{q_1} \right)+s_1 \nonumber
\\ \quad & =q_1(g_1^{-1}\cdot f)\left(\frac{t-t_1}{q_1^2},\frac{x-x_1}{q_1} \right)+q_1 s_2+s_1 \nonumber
\\ \quad & =f(t,x)+q_1 s_2+s_1 \nonumber
\\ \quad & =(\theta_1(g_1 g_2)\cdot f)(t,x). 
\end{align}
We then define $\gamma_1:G_1\times I'_1 \to I'_1$ by 
\begin{equation*}
\gamma_1(g)\cdot f(t,x)=\theta_1(g)\cdot (g\cdot f)(t,x).
\end{equation*}
Using \eqref{ThetaProduct} we obtain
\begin{align}\label{gammaction}
\gamma_1(g_1 g_2).f(t,x)&=\theta_1(g_1g_2)\cdot (g_1 g_2\cdot f)(t,x) \nonumber \\
\quad & = \theta_1(g_1) \cdot (g_1\cdot (\theta_1(g_2)\cdot(g_1^{-1}\cdot(g_1 g_2\cdot f))))(t,x) \nonumber \\
\quad & = \theta_1(g_1) \cdot (g_1\cdot (\theta_1(g_2)\cdot(g_2\cdot f)))(t,x) \nonumber \\
\quad & = \theta_1(g_1) \cdot (g_1\cdot (\gamma_1(g_2)\cdot f))(t,x) \nonumber \\
\quad & = \gamma_1(g_1) \cdot (\gamma_1(g_2)\cdot f))(t,x). 
\end{align}
Therefore $\gamma_1$ defines an action of $G_1$ on $I'_1$. Before showing that this action globalizes the local action of the one-parameter groups generated by $\mathfrak{g}$ we will write out the action explicitly. If $g_1\in G_1$ is given as in \eqref{GroupElements}, then
\begin{align}\label{NonLinearAction}
\gamma_1(g_1).f(t,x)&=\theta_1(g_1)\cdot (g_1 \cdot f)(t,x) \nonumber \\
\quad & =  (g_1 \cdot f)(t,x)+ s_1  \nonumber\\
\quad & =  q_1 f\left(\frac{t-t_1}{q_1^2},\frac{x-x_1}{q_1} \right)+ s_1 
\end{align}

\begin{proposition}
The action of $G_1$ on $I'_1$ defined by $\gamma_1$ globalizes the local action of the one parameter group actions generated by the basis elements of $\mathfrak{g}$. 
\begin{proof}
With $\xi_1, \xi_2, \xi_3,$ and $\xi_4$ as in \eqref{LAgenerators} and using \eqref{NonLinearAction} we calculate:
\begin{align*}
(\gamma_1(\exp(\epsilon \xi_1))\cdot f )(t,x)&= f(t-\epsilon,x)\\
(\gamma_1(\exp(\epsilon \xi_2))\cdot f )(t,x)&= f(t,x-\epsilon)\\
(\gamma_1(\exp(\epsilon \xi_3))\cdot f )(t,x)&= f(t,x)+\epsilon\\
(\gamma_1(\exp(\epsilon \xi_4))\cdot f )(t,x)&= e^\epsilon f(e^{-2\epsilon}t ,e^{-\epsilon}x).
\end{align*}
By taking the derivative at $\epsilon=0$ we obtain:
\begin{align*}
\frac{d}{d\epsilon}\Big|_{\epsilon=0}(\gamma_1(\exp(\epsilon \xi_1))\cdot f )(t,x)&= -f_t(t,x)\\
\frac{d}{d\epsilon}\Big|_{\epsilon=0}(\gamma_1(\exp(\epsilon \xi_2))\cdot f )(t,x)&= -f_x(t,x)\\
\frac{d}{d\epsilon}\Big|_{\epsilon=0}(\gamma_1(\exp(\epsilon \xi_3))\cdot f )(t,x)&= 1 \\
\frac{d}{d\epsilon}\Big|_{\epsilon=0}(\gamma_1(\exp(\epsilon \xi_4))\cdot f )(t,x)&= -2t f_t(t,x)-x f_x(t,x)+f(t,x).
\end{align*}
Now, it suffices to note that a vector field 
$$V=\tau(t,x)\frac{\partial}{\partial t}+\xi(t,x)\frac{\partial}{\partial x}+\phi(t,x,v(t,x))\frac{\partial}{\partial v}$$
generates a one-parameter group $G_s$ such that its action on a function $v(t,x)$ differentiated at $s=0$ gives
$$-\tau(t,x)\frac{\partial v}{\partial t}-\xi(t,x)\frac{\partial v}{\partial x}+\phi(t,x,v(t,x))$$
(c.f. \cite{Olver}). This shows that the action of $\xi_i$ corresponds to the action of $X_i$ on $I'_1$ for $1\leq i \leq 4$ and concludes the proof.
\end{proof}
\end{proposition}
%----------------------------------------------- Section 3 -------------------------------------------------------------------
\section{Case $k(v_x)=e^{v_x}$}
In this case the Lie algebra $\mathfrak{g}$ is extended by the operator $$X_5=t\frac{\partial}{\partial t}-x\frac{\partial}{\partial v}$$
(c.f. \cite{CNC}). We will denote this extended Lie algebra by $\mathfrak{g}_2$. For this case, we consider the solvable group $$G_2:=\left\{\begin{pmatrix}e^r q^2 & 0 & 0 & t \\ 0 & q & 0 & x \\ 0 & -r q& q& s \\ 0& 0 & 0 & 1\end{pmatrix} \ \Bigg| \ q, r,s,t, x  \in \R \text{ and } q>0\right\}$$
and the subgroups
$$P_2:=\left\{\begin{pmatrix} e^r q^2 & 0 & 0 & 0 \\ 0 & q &0 & 0 \\  0 & -rq &q & s \\ 0 & 0 &0 & 1 \end{pmatrix} \ \bigg| \ r, s\in \R \text{ and } q \in \R^{>0} \right\}$$
and $N_2:=N_1$.
We define the character $\chi_{2}:P_2\to \C^\times$ by $$\chi_{2}\begin{pmatrix} e^r q^2 & 0 & 0 & 0 \\ 0 & q &0 & 0 \\  0 & -rq &q & s \\ 0 & 0 &0 & 1 \end{pmatrix}=q$$
and the induced representation associated to this character would be
$$I_{2}:=\text{Ind}_{P_2}^{G_2}(\chi_{2}):=\left\{\varphi \in C^\infty(G) \ | \  \varphi(gp)=\chi_{2}(p)^{-1} \varphi(g) \text{ for all } g\in G_2 \text{ and } p\in P_2\right\}$$
with the $G_2$ action of left translation. Restriction of domain allows us to define 
$$I'_{2}:=\left\{ f\in C^\infty(\R^2) \ | \ f(t,x)=\varphi(n_{t,x}) \text{ for some } \varphi \in I_2 \text{ and for all } t,x\in \R^2 \right\}.$$
For the rest of this section let
\begin{align*}
g_1&:=\begin{pmatrix}e^{r_1} q_1^2 & 0 & 0 & t_1 \\ 0 & q_1 &0 & x_1 \\  0 & -r_1 q_1 &q_1 & s_1 \\ 0 & 0 &0 & 1 \end{pmatrix}, & g_2&:=\begin{pmatrix}e^{r_2} q_2^2 & 0 & 0 & t_2 \\ 0 & q_2 &0 & x_2 \\  0 & -r_2 q_2 &q_2 & s_2 \\ 0 & 0 &0 & 1 \end{pmatrix}.
\end{align*}
Notice that
$$
\begin{pmatrix}e^r q^2 & 0 & 0 & t \\ 0 & q &0 & x \\  0 & -r q &q & s \\ 0 & 0 &0 & 1 \end{pmatrix}=\begin{pmatrix} 1 & 0 & 0 &t  \\ 0 & 1&0 & x \\  0 & 0 &1 & 0 \\ 0 & 0 &0 & 1 \end{pmatrix}\begin{pmatrix} e^r q^2 & 0 & 0 & 0 \\ 0 & q &0 & 0 \\  0 & -r q &q & s \\ 0 & 0 &0 & 1 \end{pmatrix} .
$$
Then for $\varphi\in I_2$ we have, $$\varphi\begin{pmatrix}e^r q^2 & 0 & 0 & t \\ 0 & q &0 & x \\  0 & -r q &q & s \\ 0 & 0 &0 & 1 \end{pmatrix}=q^{-1}f(t,x).$$
As before, we can use this to show that $I_2\cong I'_2 \cong C^\infty(\R^2)$ as vector spaces. In order for the restriction map to become intertwining, we define the action of $G_2$ on $I'_2$ by
\begin{equation*}\label{ActionOfG2}g_1\cdot f(t,x)=q_1 f\left(\frac{t-t_1}{e^{r_1} q_1^2},\frac{x-x_1}{q_1} \right).\end{equation*}
To globalize the action of the one parameter groups generated by $\mathfrak g$ we define the map $\theta_2:G_2\times I'_2\to I'_2$ by
\begin{equation*}(\theta_2(g_1)\cdot f)(t,x)=f(t,x)+s_1 -r_1(x-x_1). \end{equation*}
As before, it is easy to see that $\theta_2$ does not define an action of $G_2$ on $I'_2$. However, we can define $\gamma_2:G_2\times I'_2\to I'_2$ by 
\begin{equation}\label{gamma2}(\gamma_2(g)\cdot f)(t,x)=(\theta_2(g)\cdot(g\cdot f))(t,x). \end{equation}
\begin{proposition}
 $\gamma_2:G_2\times I'_2\to I'_2$ defined as in \eqref{gamma2} defines an action of $G_2$ on $I'_2$.
\begin{proof}
Notice that $$g_1g_2=\begin{pmatrix}e^{r_1+r_2} (q_1q_2)^2 & 0 & 0 & t_1+e^{r_1}q_1 t_2 \\ 0 & q_1q_2 &0 & x_1+q_1x_2 \\  0 & -(r_1+r_2) q_1q_2 &q_1 q_2 & s_1+q_1(s_2-r_1x_2) \\ 0 & 0 &0 & 1 \end{pmatrix}.$$
Then, 
\begin{align*}(\theta_2(g_1 g_2)\cdot f)(t,x)&=f(t,x)+s_1+q_1(s_2-r_1x_2)-(r_1+r_2)(x-( x_1+q_1x_2))
\\ \quad & =f(t,x)+s_1+q_1 s_2-r_1(x-x_1)+r_2(x- x_1-q_1x_2).
\end{align*}
On the other hand,
\begin{align*}
\theta_2(g_1) \cdot (g_1\cdot (\theta_2(g_2)\cdot(&g_1^{-1}\cdot f)))(t,x) =g_1\cdot (\theta_2(g_2)\cdot(g_1^{-1}\cdot f))(t,x)+s_1-r_1(x-x_1) \nonumber
\\ \quad &=q_1(\theta_2(g_2)\cdot(g_1^{-1}\cdot f))\left(\frac{t-t_1}{q_1^2},\frac{x-x_1}{q_1} \right)+s_1-r_1(x-x_1) \nonumber
\\ \quad & =f(t,x)+s_1-r_1(x-x_1) \nonumber+ q_1(s_2-r_2(x-(x-x_1)/q_1)) \nonumber
\\ \quad & =(\theta_2(g_1 g_2)\cdot f)(t,x). 
\end{align*}
Now, an identical calculation as in \eqref{gammaction} replacing the subindex $1$ for $2$ gives the desired result.
\end{proof}
\end{proposition}
The global action $\gamma_2$ of $G_2$ on $I'_2$ is given by
\begin{equation}\label{Nonlinear2}(\gamma_2(g_1)\cdot f)(t,x)=q_1 f\left(\frac{t-t_1}{e^{r_1} q_1^2},\frac{x-x_1}{q_1} \right)+s_1-r_1(x-x_1) \end{equation}
\begin{proposition}
The action of $G_2$ on $I'_2$ defined by $\gamma_2$ globalizes the local action of the one parameter group actions generated by the basis elements of $\mathfrak{g_2}$. 
\begin{proof}
Let $\xi_1, \xi_2, \xi_3,$ and $\xi_4$ be defined as in \eqref{LAgenerators} and let $$\xi_5:= \begin{pmatrix}1 & 0 & 0 & 0 \\ 0 & 0 &0 & 0 \\  0 & -1  &0 & 0 \\ 0 & 0 &0 & 0 \end{pmatrix}.$$

Then, using Equation \eqref{Nonlinear2} we obtain:
\begin{align*}
(\gamma_2(\exp(\epsilon \xi_1))\cdot f )(t,x)&= f(t-\epsilon,x)\\
(\gamma_2(\exp(\epsilon \xi_2))\cdot f )(t,x)&= f(t,x-\epsilon)\\
(\gamma_2(\exp(\epsilon \xi_3))\cdot f )(t,x)&= f(t,x)+\epsilon\\
(\gamma_2(\exp(\epsilon \xi_4))\cdot f )(t,x)&= e^\epsilon f(e^{-2\epsilon}t ,e^{-\epsilon}x) \\
(\gamma_2(\exp(\epsilon \xi_5))\cdot f )(t,x)&= f(e^{-\epsilon}t ,x)-\epsilon x.
\end{align*}
By taking the derivative at $\epsilon=0$ we obtain:
\begin{align*}
\frac{d}{d\epsilon}\Big|_{\epsilon=0}(\gamma_2(\exp(\epsilon \xi_1))\cdot f )(t,x)&= -f_t(t,x)\\
\frac{d}{d\epsilon}\Big|_{\epsilon=0}(\gamma_2(\exp(\epsilon \xi_2))\cdot f )(t,x)&= -f_x(t,x)\\
\frac{d}{d\epsilon}\Big|_{\epsilon=0}(\gamma_2(\exp(\epsilon \xi_3))\cdot f )(t,x)&= 1 \\
\frac{d}{d\epsilon}\Big|_{\epsilon=0}(\gamma_2(\exp(\epsilon \xi_4))\cdot f )(t,x)&= -2t f_t(t,x)-x f_x(t,x)+f(t,x) \\
\frac{d}{d\epsilon}\Big|_{\epsilon=0}(\gamma_2(\exp(\epsilon \xi_5))\cdot f )(t,x)&= -f_t(t ,x)- x.
\end{align*}
This finishes the proof, as it shows that the one-parameter group generated by $X_i$ is the same as the one generated by $\xi_i$ for each $1\leq i \leq 5$.
\end{proof}
\end{proposition}

%-------------------------------------------------------------- Section 4 ---------------------------------------------------
\section{Case $k(v_x)=(v_x)^{n}, n\geq -1, n\neq 0$}

In this case the Lie algebra $\mathfrak{g}$ is extended by the operator $$X_6=nt\frac{\partial}{\partial t}-v\frac{\partial}{\partial v}$$
(c.f. \cite{CNC}). Denote this Lie algebra by $\mathfrak{g}_3$. Now we will consider the solvable group
 $$G_3:=\left\{\begin{pmatrix}r^n q^2 & 0 & 0 & t \\ 0 & q & 0 & x \\ 0 & 0 & q r^{-1}& s \\ 0& 0 & 0 & 1\end{pmatrix} \ \Bigg| \ q, r,s,t, x  \in \R \text{ and }  q,r >0\right\}$$
with subgroups 
$$P_3:=\left\{\begin{pmatrix}r^n q^2 & 0 & 0 & 0 \\ 0 & q & 0 & 0 \\ 0 & 0 & q r^{-1}& s \\ 0& 0 & 0 & 1\end{pmatrix}\ \Bigg| \ q, r,s  \in \R \text{ and } q,r>0\right\}$$ and 
$N_3:=N_1$.
As in the previous cases, we define the character $\chi_{3}:P_3\to \C^\times$ by $$\chi_{3}\begin{pmatrix}r^n q^2 & 0 & 0 & 0 \\ 0 & q & 0 & 0 \\ 0 & 0 & q r^{-1}& s \\ 0& 0 & 0 & 1\end{pmatrix}=q r^{-1}$$
and use it to induce a representation of $G_3$ in the standard way
$$I_{3}:=\text{Ind}_{P_3}^{G_3}(\chi_{3}):=\left\{\varphi \in C^\infty(G) \ | \  \varphi(gp)=\chi_{3}(p)^{-1} \varphi(g) \text{ for all } g\in G_3 \text{ and } p\in P_3\right\}.$$
Using the fact that 
$$\begin{pmatrix}r^n q^2 & 0 & 0 & t \\ 0 & q & 0 & x \\ 0 & 0 & q r^{-1}& s \\ 0& 0 & 0 & 1\end{pmatrix}=\begin{pmatrix} 1 & 0 & 0 &t  \\ 0 & 1&0 & x \\  0 & 0 &1 & 0 \\ 0 & 0 &0 & 1 \end{pmatrix}\begin{pmatrix}r^n q^2 & 0 & 0 & 0 \\ 0 & q & 0 & 0 \\ 0 & 0 & q r^{-1}& s \\ 0& 0 & 0 & 1\end{pmatrix},$$
it is easy to show that $I_3$ is isomorphic as a vector space to
$$I'_{3}:=\left\{ f\in C^\infty(\R^2) \ | \ f(t,x)=\varphi(n_{t,x}) \text{ for some } \varphi \in I_3 \text{ and for all } t,x\in \R^2 \right\}$$
and that both are isomorphic to $ C^\infty(\R^2)$. 
For the rest of the section let 
\begin{align*}\label{GroupElements4} g_1=&\begin{pmatrix} r_1^n q_1^2 & 0 & 0 & t_1 \\ 0 &  q_1 &0 & x_1 \\  0 & 0 &r_1^{-1}q_1 & s_1 \\ 0 & 0 &0 & 1 \end{pmatrix}. \end{align*} 
We give $I'_3$ the structure of a $G_3$-representation by defining the action of $G_3$ by
$$g_1\cdot f(t,x)=r_1^{-1} q_1 f\left(\frac{t_2-t_1}{r_1^n q_1^2},\frac{x_2-x_1}{q_1} \right)$$
so that $I_3\cong I'_3$ as $G_3$-representations. The maps to obtain the globalization for this particular case are very similar to the ones in Section \ref{genericcase}, thus most of the details will be omitted. Suffices to say that $\theta_3: G_3\times I'_3 \to I'_3$ is given by $\theta_3(g_1)\cdot f(t,x)=f(t,x)+s_1$ and that the map $\gamma_3: G_3\times I'_3 \to I'_3$ given by $\gamma_3(g)\cdot f(t,x)=(\theta_3(g)\cdot (g\cdot f))(t,x)$ defines an action of $G_3$ on $I'_3$. In particular, this action is defined by:
\begin{equation}\label{ActionGamma3} \gamma_3(g_1)\cdot f(t,x)=r_1^{-1} q_1 f\left(\frac{t_2-t_1}{r_1^n q_1^2},\frac{x_2-x_1}{q_1} \right)+s_1. \end{equation}

\begin{proposition}
The action $\gamma_3: G_3\times I'_3 \to I'_3$ globalizes the local action of the one parameter group actions generated by the basis elements of $\mathfrak{g_3}$. 
\begin{proof}
With $\xi_1, \xi_2, \xi_3,$ and $\xi_4$ defined as in \eqref{LAgenerators} and $\xi_6$ given by $$\xi_6:= \begin{pmatrix}n & 0 & 0 & 0 \\ 0 & 0 &0 & 0 \\  0 & 0  &-1 & 0 \\ 0 & 0 &0 & 0 \end{pmatrix}.$$
we can use Equation \eqref{ActionGamma3} to obtain:
\begin{align*}
(\gamma_3(\exp(\epsilon \xi_1))\cdot f )(t,x)&= f(t-\epsilon,x)\\
(\gamma_3(\exp(\epsilon \xi_2))\cdot f )(t,x)&= f(t,x-\epsilon)\\
(\gamma_3(\exp(\epsilon \xi_3))\cdot f )(t,x)&= f(t,x)+\epsilon\\
(\gamma_3(\exp(\epsilon \xi_4))\cdot f )(t,x)&= e^\epsilon f(e^{-2\epsilon}t ,e^{-\epsilon}x) \\
(\gamma_3(\exp(\epsilon \xi_6))\cdot f )(t,x)&= e^{-\epsilon}f(e^{-n \epsilon}t ,x)
\end{align*}
and taking the derivative at $\epsilon=0$ we obtain:
\begin{align*}
\frac{d}{d\epsilon}\Big|_{\epsilon=0}(\gamma_3(\exp(\epsilon \xi_1))\cdot f )(t,x)&= -f_t(t,x)\\
\frac{d}{d\epsilon}\Big|_{\epsilon=0}(\gamma_3(\exp(\epsilon \xi_2))\cdot f )(t,x)&= -f_x(t,x)\\
\frac{d}{d\epsilon}\Big|_{\epsilon=0}(\gamma_3(\exp(\epsilon \xi_3))\cdot f )(t,x)&= 1 \\
\frac{d}{d\epsilon}\Big|_{\epsilon=0}(\gamma_3(\exp(\epsilon \xi_4))\cdot f )(t,x)&= -2t f_t(t,x)-x f_x(t,x)+f(t,x) \\
\frac{d}{d\epsilon}\Big|_{\epsilon=0}(\gamma_3(\exp(\epsilon \xi_6))\cdot f )(t,x)&= -f(t ,x)- n f_t(t,x).
\end{align*}
Now, the proof is completed as in the previous cases.
\end{proof}
\end{proposition}

\section{Case $k(v_x)=\frac{e^{n \arctan(v_x)}}{1+v_x^2}$}

In this case the Lie algebra $\mathfrak{g}$ is extended by the operator $$X_7=nt\frac{\partial}{\partial t}-x\frac{\partial}{\partial v}+v\frac{\partial}{\partial x}$$
(c.f. \cite{CNC}). In this particular case a globalization for the action of $X_7$ on a non-trivial space of functions cannot be constructed. To see this, suppose that $G$ is a Lie group that globalizes the local action generated by $X_7$ and let $g\in G$. Consider a function $f $ defined on a subset $\Omega\subset \R^2$ that contains an open set and consider its graph $\Gamma_f:=\{(t, x,f(t, x)) \ | \ (t,x)\in \Omega\}$. If $\Gamma_f\subset \R^2\times \R$ is invariant under the action of $g$, then
$$g\cdot \Gamma_f=\{ (\tilde{t},\tilde{x},\tilde{v})=g\cdot(t,x,v) \  | \ v=f(t,x)\}.$$  
In particular, by exponentiating the standard action of $X_7$ on a $\R^2\times \R$,  we obtain
\begin{equation}\label{LocalAction}\exp(\epsilon X_7)\cdot \Gamma_f=\{ (e^{n\epsilon}t,x\cos \epsilon+v\sin \epsilon, v\cos\epsilon-x\sin \epsilon) \  | \ v=f(t,x)\}.\end{equation}
Since the action includes a rotation by an angle $\epsilon$ on the $x,v$-coordinates and a dilation in the $t$-coordinate, unless the function $f$ does not depend on $x$ we cannot guarantee  that $\exp(\epsilon X_7)\cdot \Gamma_f$ is the graph of a single-valued function for all $\epsilon \in \R$ for a general function $f$. 

However, this does not mean that $X_7$ does not generate a local action on a function $f$. To illustrate, let $$f(x)=ax+b.$$
Then, $f$ is a solution of \eqref{NLFE}, its graph is $\Gamma_f=\{(t,x,ax+b) | t,x\in \R\}$, and from \eqref{LocalAction} we have $$\exp(\epsilon X_7)\cdot(t,x,ax+b)=(t,x\cos \epsilon +(ax+b)\sin\epsilon, (ax+b)\cos\epsilon -x \sin \epsilon).$$
From this, we can obtain that

$$(\exp(\epsilon X_7)\cdot f)(\tilde{t},\tilde{x})=\frac{a\cos\epsilon-\sin\epsilon}{a\sin \epsilon + \cos \epsilon}\tilde{x}+\frac{b}{a\sin \epsilon + \cos \epsilon}$$
which is a linear function, provided  $a\neq -\cot \epsilon$. So $X_7$ generates a local action on $f$, but it cannot be globalized.

%\bibliographystyle{plain}
%\bibliography{Bib}
\end{document}